\documentstyle[12pt]{article}
\newlength{\abstractwidth}
\renewcommand{\thefootnote}{\fnsymbol{footnote}}
\renewcommand{\thanks}[1]{\footnote{#1}} 
\newcommand{\starttext}{ \setcounter{footnote}{0}
\renewcommand{\thefootnote}{\arabic{footnote}}}

\newcommand{\be}{\begin{equation}}
\newcommand{\bea}{\begin{eqnarray}}
\newcommand{\eea}{\end{eqnarray}} \newcommand{\ee}{\end{equation}}
 
 \def\ba{\begin{eqnarray}}
\def\ea{\end{eqnarray}}

\newcommand{\db}{\overline{\partial}}

\def\o{\omega}
\def\Re{{\rm Re}}

\def\tr{{\rm tr}}
\def\det{{\rm det}}

\def\log{\,{\rm log}\,}

\def\o{\omega}

\def\o{\omega}

\def\vp{\varphi}

\def\ge{\geq}
\def\le{\leq}

\def\R{{\bf R}}

\def\p{\partial}

\def\ma{{Monge-Amp\`ere\ }}

\def\[{{\bf [}}
\def\]{{\bf ]}}

\def\ddbar{i\p\bar\p}

\def\mathbb{\bf}

\def\eqref{\ref}

\newcommand{\neweqref}[1]{(\ref{#1})}

\begin{document}
\starttext \baselineskip=18pt \setcounter{footnote}{0}
\newtheorem{theorem}{Theorem}
\newtheorem{lemma}{Lemma}
\newtheorem{corollary}{Corollary}
\newtheorem{definition}{Definition}
\newtheorem{conjecture}{Conjecture}
\newtheorem{proposition}{Proposition}


\begin{center}
{\Large \bf On uniform estimates for $(n-1)-$form fully nonlinear partial differential equations on compact Hermitian manifolds}

\medskip
\centerline{Nikita Klemyatin, Shuang Liang and Chuwen Wang}

\medskip

\begin{abstract}

{\footnotesize We obtain a priori $L^\infty$ estimate for a general class of $(n-1)-$form fully nonlinear partial differential equations on compact Hermitian manifolds. Our method relies on the local version of comparison with auxiliary Monge-Amp\`ere equations, developed earlier by B. Guo and D. H. Phong. The key is to find the appropriate elliptic operator such that the maximum principle applies.}

\end{abstract}

\end{center}

\baselineskip=15pt
\setcounter{equation}{0}
\setcounter{footnote}{0}

\section{Introduction}
\setcounter{equation}{0}

A priori estimates have always been fundamental in the study of partial differential equations. Among them, $L^\infty$ estimates are especially important. As the most prominent example, $L^\infty$ estimate for the \ma equation obtained through Moser iteration was crucial in Yau's solution to Calabi conjecture \cite{Y}. Kolodziej \cite{K} improved it to a sharp version using pluripotential theory. A recent development initiated by Guo, Phong and Tong \cite{GPT} provided a PDE-based proof to Kolodziej's result using comparison with solutions to auxiliary \ma equations. Moreover, the method extends to a wide class of non-linear equations satisfying a structural condition, which has been shown by Harvey and Lawson \cite{HL3} to be quite large. This flexible method can also be applied to obtain stability estimates for Monge-Amp\`ere and Hessian equations \cite{GPTa}; $L^\infty$ estimates for Monge-Amp\`ere equations on nef classes rather than just K\"ahler classes \cite{GPTW}; sharp modulus of continuity for non-H\"older solutions \cite{GPTW1}; extensions to parabolic equations \cite{CC2}; Regularization of $m-$subharmonic functions and h\"older continuity \cite{CX}; lower bounds for the Green's function \cite{GPS}; uniform entropy estimates \cite{GPa}; and diameter estimates and convergence theorems in K\"ahler geometry not requiring bounds on the Ricci curvature \cite{GPSS,GS}. Sroka \cite{S} applied the same method to obtain a sharp uniform bound for the quaternionic Monge-Amp\'ere equation on hyperhermitian manifolds. For an exposition of these topics, we refer the readers to the recent survey paper by Guo and Phong \cite{GPb}.

Recently, Guo and Phong \cite{GP} developed a local version of the comparison method to extend the $L^\infty$ bound to equations on Hermitian manifolds. It also applies to $(n-1)-$form \ma equations \cite{GP}, which was first solved by Tosatti-Weinkove \cite{TW1} with regard to Calabi-Yau equation for Gauduchon metrics. A natural question is whether the comparison method in \cite{GP} applies to more general $(n-1)-$form fully nonlinear equations other than \ma equations. In this paper, we generalize the $L^\infty$ estimate to a wide class of $(n-1)-$form fully nonlinear equations satisfying structure conditions proposed in \cite{GPT}. The key is to find the appropriate elliptic operator similar to the one defined in \cite[p.17]{TW1} and \cite[p.20]{GP}, so that the maximum principle applies to the test function.

We now describe our results in detail.

Consider a positive function $f: \Gamma \subset {\mathbb R}^n\to {\mathbf R}_+$ such that 

 (1) the domain $\Gamma\subset {\mathbb R}^n$ is an open symmetric cone satisfying 
\begin{equation}\label{eqn:cone}
\Gamma_n\subset \Gamma \subset \Gamma_1;
\end{equation}
Here $ \Gamma_k:=\{\lambda \in \R^n :\sigma_j(\lambda)>0,\; 1\leq j\leq k\} $ is the Garding cones, where $\sigma_j(\lambda)$ is the $j$-th symmetric polynomial in $\lambda$.

(2) $f(\lambda)$ is symmetric in $\lambda = (\lambda_1,\ldots, \lambda_n)\in \Gamma$ and it satisfies 
\bea\label{eqn:eulerineq}
\sum_j \frac{\partial f}{\partial \lambda_j}\lambda_j \leq \beta f
\eea
for some positive constant $\beta>0$;

 (3) $\frac{\partial f}{\partial \lambda_j}>0$ for all $\lambda\in \Gamma$, $j=1,\ldots, n$;

 (4) There exists $\gamma>0$ such that 
\begin{equation}\label{eqn:structure}
\prod_{j=1}^n \frac{\partial f}{\partial \lambda_j}\ge \gamma.
\end{equation}
for all $\lambda\in \Gamma$.

Remark that (\ref{eqn:eulerineq}) slightly relaxes the conditions in \cite{GPT,GPTa,GPTW,GPTW1,GPS,GP}, which require $f$ to be homogeneous of degree one.









Suppose $(X,\omega)$ is a compact Hermitian manifold without boundary and $\omega_h$ is another Hermitian metric on $X$. For any $\varphi \in C^2(X)$, we set 
\bea
\label{eqn:unknownmetric}
\tilde \o = \omega_h + \frac{1}{n-1} ( (\Delta_\omega \varphi) \omega - \ddbar \varphi) + \chi[\vp],
\eea
where $$\Delta_\omega \varphi = n\frac{\ddbar \varphi \wedge \omega^{n-1}}{\omega^n}$$ is the complex Laplacian of $\varphi$ with respect to $\omega$ and $\chi[\vp]=x(z,\vp,\partial \vp, \db \vp)$ is the gradient term which may depend on the $\varphi$ and its first order derivatives.

Let $h_\varphi: TX\to TX$ be the relative endomorphism from $\omega$ to $\tilde{\o}$. In local coordinates, $(h_\varphi)^j{}_i= g^{j\bar k} \tilde g_{\bar k i}$, where $(\tilde g_{\bar k i})$ denotes the components of $\tilde \o$ and $(g^{i\bar j})$ denotes the inverse of $\o$. Let $\lambda[h_\varphi]$ be the (unordered) vector of eigenvalues of $h_\varphi$.

For a smooth function $F$ on $X$, we consider the fully nonlinear partial differential equation
\begin{equation}\label{eqn:main}
f(\lambda[h_\varphi]) = e^{F},
\mbox{ and } \lambda[h_\varphi]\in \Gamma. 
\end{equation}
with $\sup_X \varphi = 0$ and $\lambda[h_\varphi] \in \Gamma$ on $X$. 

We remark that our definition of relative endomorphism $h_\varphi := \o ^{-1} \tilde \o$ is different from those in \cite{GPT,GPTa,GPTW,GPTW1,GPS,GP} due to a different form of unknown metric \neweqref{eqn:unknownmetric}. So the condition $\lambda \in \Gamma$ means $\lambda[\o^{-1}\cdot \tilde \o] \in \Gamma $, instead of $\lambda[\o^{-1}\cdot (\o + i\partial \bar \partial \varphi)]\in \Gamma$, which imposes different condition on $\varphi$ as in \cite{GPT,GPTa,GPTW,GPTW1,GPS,GP}. For example, when $\chi = 0, \, \Gamma = \Gamma_n$, this requires $\varphi$ being $(n-1)-$plurisubharmonic in the sense of \cite{HLnew,HLnew2}, instead of being plurisubharmonic.

The following two special cases is of particular interests. When $f(\lambda[h_\varphi]) = \Big(\frac{\tilde \o^n}{\omega^n} \Big)^{1/n}$, $\chi = 0$ and $\Gamma = \Gamma_n$, the equation is the $(n-1)-$form Monge-Amp\`ere equation considered by Tosatti-Weinkove \cite{TW1}. They showed  when $\omega$ is K\"ahler, the equation is solvable if $F$ is modified by a suitable additive constant, in which the $C^0$ estimate is an essential step of their proof. When $f(\lambda[h_\varphi]) = \Big(\frac{\tilde \o^n}{\omega^n} \Big)^{1/n}$, $\chi[\vp]=* \Re \{\sqrt{-1}\partial \vp \wedge \bar{\partial} \omega^{n-2}\}$ and $\Gamma = \Gamma_n$, the equation is considered by Sz\'ekelyhidi-Tosatti-Weinkove \cite{STW} in their celebrated solution to the Gauduchon conjecture. 

Later in Guo-Phong \cite{GP}, {\em a priori} $C^0$ estimate of $\varphi$ can be obtained for the $(n-1)-$form Monge-Amp\`ere equation even when $\o$ is not K\"ahler. Here we generalize their results to arbitrary equations $f$ satistying the conditions (1-4) above with the assumption the gradient term $\chi$ being semi-positive.

\begin{theorem}
\label{local}\label{thm:local}
Assume $\chi \geq 0$, let $\varphi$ be a $C^2$ solution on a compact Hermitian manifold $(X,\o)$ of the equation
\bea
\label{eqn:f}
f(\lambda[h_\varphi])=e^F
\eea
where the operator $f(\lambda)$ satisfies the conditions (1-4) spelled out in Section \S 1. 

Fix any $p>n$, we have
\bea
||\varphi||_{L^\infty(X)} \leq C
\eea
where $C$ is a constant depending only on $X,\o,\o_h, n,p$, and $\|e^{F}\|_{L^1(\log L)^p}$. Here the $L^1(\log L)^p$ norms are with respect to the volume form $\o^n$.
\end{theorem}

Recently, Guo-Phong manage to prove $C^0$ estimates for the gradient term $\chi$ of a specific form without being semi-positive. Their method uses real auxiliary Monge-Amp\`ere equation, we refer the interested readers to the most recent version of \cite{GP}.

\section{The local estimate and proof of Theorem 1}\label{section 2 new}
\setcounter{equation}{0}

We denote by $G^{i\bar j} = \frac{\partial \log f(\lambda[h])}{\partial h_{ij}} = \frac{1}{f} \frac{\partial f(\lambda[h])}{\partial h_{ij}}$ the coefficients of the linearization of the operator $\log f(\lambda[h])$.

It follows from the structure conditions of $f$ that $G^{i\bar j}$ is positive definite at $h_{\varphi}$ and
\bea
\label{eqn:structure}
\det ( G^{i\bar j}) = \frac {1}{f^n} \det (\frac{\partial f(\lambda[h])}{\partial h_{ij}}) \ge \frac{\gamma}{f(\lambda)^n}.
\eea

Fix any point $z_0 \in X$, by simple linear algebra there is a smooth local frame $\{e^i\}$  for the holomorphic cotangent bundle ${T^*}^{1,0}X$ in a neighborhood of $z_0$, such that $\o = \sqrt{-1} \sum_j e^j \wedge e^{\bar{j}}$ and $\tilde \o =  \sqrt{-1} \sum_j \lambda_j  e^j \wedge e^{\bar{j}}$. In particular, $(h_\varphi)_{\bar j i} = \lambda_j \delta_{ij}$ and $G^{i\bar j} = \frac{1}{f} \frac{\partial f}{\partial \lambda_j}\delta_{ij}$ in this frame. Therefore, we have in the neighborhood of $z_0$
\bea
\label{eqn:homogeneity}
\tr_G \tilde\omega = \sum_{i,j} G^{i\bar j} \tilde g_{\bar j i} = \sum_j \frac{1}{f} \frac{\partial f}{\partial \lambda_j}\lambda_j \leq \beta.  
\eea
The last equality follows from (\ref{eqn:eulerineq}). As $z_0$ is arbitrary, we know that $\tr_G \tilde \o \leq C$ holds everywhere on $X$.

Define the following tensor $\Theta^{i\bar j}$ by 
$$\Theta^{i\bar j} = \frac{1}{n-1} \Big( (G^{k\bar l} g_{\bar l k}) g^{i\bar j}  -  G^{i\bar j}   \Big).$$

\noindent\textbf{Remark.} This tensor is similar to the tensor $\Theta^{i\bar j}$ defined in \cite[p.17]{TW1} and \cite[p.20]{GP}, except that we put $G^{i\bar j}$ in place of $\tilde g^{i\bar j}$ in order to deal with more general fully nonlinear equations.

We summarize the key properties of $\Theta^{i\bar j}$ in the following lemma.
\begin{lemma}
\label{lm:theta}
The tensor $\Theta^{i\bar j}$ then satisfies the following: for $\lambda\in \Gamma$,

{\rm (a)} $\Theta^{i\bar j}\varphi_{i\bar j} \leq \beta - G^{i\bar j}(g_h)_{\bar j i} -  G^{i\bar j}\chi_{\bar j i} $ and $G^{i\bar j}(g_h)_{\bar j i}, G^{i\bar j}\chi_{\bar j i} \geq 0$;

{\rm (b)} $\Theta^{i\bar j}$ is positive definite, and $\det(\Theta^{i\bar j}) \ge \frac{\gamma}{f(\lambda)^n}$.

\end{lemma}

\noindent{\em Proof. } 
We have
\bea 
\beta &\geq& \tr_G \tilde\omega \nonumber \\
      &=& G^{i\bar j} \Big( (g_h)_{\bar j i} + \frac{1}{n-1}(g^{k\bar l}\varphi_{\bar k l} g_{\bar j i} - \varphi_{\bar j i}) + \chi_{\bar j i} \Big)  \nonumber\\
      &=& G^{i\bar j}(g_h)_{\bar j i} + G^{i\bar j}\chi_{\bar j i}  + \frac{1}{n-1} \Big(  (G^{k\bar l} g_{\bar l k}) g^{i\bar j} - G^{\bar j i} \Big) \varphi_{\bar j i}  \nonumber \\
      &=& G^{i\bar j}(g_h)_{\bar j i} + G^{i\bar j}\chi_{\bar j i}  + \Theta^{i\bar j} \varphi_{\bar j i} \nonumber
\eea
where the first equality follows from the equation \neweqref{eqn:homogeneity}, the second equality follows from definition of $\tilde \o$ and the last equality follows from the definition of $\Theta$. Since $G^{i\bar j}$ is positive definite, we know $G^{i\bar j}(g_h)_{\bar j i} \geq 0$. This proves (a).

We can choose holomorphic coordinates at a given point $p\in X$ such that $g_{i\bar j}|_{p} = \delta_{ij} $ and $ G_{i\bar j}|_{p} = \mu_i \delta_{ij}$. Note that $\mu_i>0$ for each $i = 1,\cdots, n$ since $G^{i\bar j}$ is positive definite. Then we have
$$\Theta^{i\bar j} |_{p} = \frac{1}{n-1} (\sum_{k\neq i} \mu_k) \delta_{ij}$$ which is clearly positive definite. Moreover, 
$$\det (\Theta^{i\bar j})|_{x_0} = \frac{1}{(n-1)^n} \prod_{i=1}^n (\sum_{k\neq i} \mu_k) \ge \prod_{i=1}^n \Big( \prod_{k\neq i} \mu_k  \Big)^{1/(n-1)} = \prod_{i=1}^n \mu_i = \det (G^{i\bar j}) \ge \frac{\gamma}{f(\lambda)^n},$$
where the middle inequality follows from the arithmetic-geometric (AG) inequality and the last inequality follows from the equation \neweqref{eqn:structure}. This proves (b). Q.E.D.

Let $x_0$ be the point where $\varphi$ attains its minimum. Without loss of generality, we assume $-\varphi(x_0)\ge 2$. We then fix a local holomorphic coordinate $z$ centered at $x_0$ such that
\bea\label{eqn:normalization}
{1\over 2}i\p\bar\p |z|^2\leq\o\leq 2 i\p\bar\p |z|^2 \ {\rm in} \ B(x_0,2r_0)
\eea
where $B(x_0,2r_0)$ is the Euclidian ball of radius $2r_0$ in this coordinate. We denote $\Omega:= B(x_0,2r_0)$ for simplicity.

Choose a small constant $\epsilon'>0$ such that 
\begin{equation}\label{eqn:epsilon small}
\omega_h \ge \frac{2\epsilon'}{n-1}(\tr_\omega \omega_h) \omega \ {\rm in} \ \Omega.
\end{equation}

Define $s_0 := 4 \epsilon' r_0^2$. Then for any $s\in (0,s_0)$, we consider the following comparison function $u_s: \Omega \to \R$ given by
\begin{equation}\label{eqn:u s 1}
u_s(z) : = \varphi(z) - \varphi(x_0) + \epsilon' |z|^2 - s,
\end{equation}
Let
\bea
\Omega_s=\{z\in \Omega; u_s(z)<0\}.
\eea
be the sublevel set of $u_s$.
It's easy to see from the definition that $u_s$ is positive on $\partial \Omega$, so sub-level sets $\Omega_s$
is relatively compact in $\Omega$.

Set
\bea
A_s=\int_{\Omega_s}(-u_s)e^{nF}\o^n.
\eea

To make the right hand side of our auxiliary \ma equation smooth, we choose the following sequence of smooth positive functions $\tau_k: {\mathbb R}\to{{\mathbb R}}_+$ such that 
\begin{equation}\label{eqn:tau k}\tau_k(x) = x + \frac 1 k,\quad \mbox{ when }x\ge 0,
\end{equation}
 and $$ \tau_k(x) = \frac 1 {2k},\quad  \mbox { when } x\le -\frac {1}{k},$$ and $\tau_k(x)$ lies between $1/2k$ and $1/k$ for $x\in [-1/k, 0]$. Clearly $\tau_k$ converge pointwise to $\tau_\infty(x) = x \cdot \chi_{{\mathbb R}_+}(x)$ as $k\to\infty$, where $\chi_{{\mathbb R}_+}$ denotes the characteristic function of ${\mathbb R}_+$.

The auxiliary Monge-Amp\`ere equation we consider is the following
\bea
\label{eqn:psi}
(i\p\bar\p \psi_{s,k})^n={\tau_k(-u_s)\over A_{s,k}}e^{nF}\o^n\ {\rm in}\ \Omega, 
\quad 
\psi_{s,k}=0\ {\rm on}\ \p\Omega
\eea
with $i\p\bar\p\psi_{s,k}\geq 0$, and $A_{s,k}$ is defined by
\bea
A_{s,k}=\int_\Omega \tau_k(-u_s) e^{nF} \omega^n.
\eea
By Caffarelli-Kohn-Nirenberg-Spruck \cite{CKNS}, this Dirichlet problem admits a unique solution $\psi_{s,k}$ which is of class $C^\infty(\bar\Omega)$, with $\psi_{s,k} \leq 0$. By the definition of $A_{s,k}$, we have $A_{s,k}\to A_s$ as $k\to\infty$, and
\bea
\int_\Omega (i\p\bar\p\psi_{s,k})^n=1.
\eea

Now we are ready to establish the following key comparison lemma

\begin{lemma}
\label{lm:comparison}
Let $u_s$ be a $C^2$ solution of the fully non-linear equation (\eqref{eqn:f}) and $\psi_{s,k}$ be the solutions of the complex Monge-Amp\`ere equation (\eqref{eqn:psi}) as defined above. Then we have
\bea
\label{comparison}
-u_s\leq \varepsilon (-\psi_{s,k})^{n\over n+1}\ {\rm on}\ \bar\Omega
\eea
where $\varepsilon$ is the constant defined by $\varepsilon^{n+1}=A_{s,k}\gamma^{-1}\beta^{n+1}{(n+1)^n\over n^{2n}}$.
\end{lemma}

\medskip
\noindent
{\it Proof.} We show that the function
\bea\label{eqn:test function}
\Phi=  -\varepsilon (-\psi_{s,k})^{n\over n+1}-u_s
\eea
is always $\leq 0$ on $\bar\Omega$. Let $x_{\max}\in \overline{\Omega}$ be a maximum point of $\Phi$. If $x_{\max}\in\overline{ \Omega}\backslash \Omega_s$, clearly $\Phi(x_{\max})\le 0$ by the definition of $\Omega_s$ and the fact that $\psi_{s,k}< 0 $ in $\Omega$. If $x_{\max}\in \Omega_s$, then we have $\ddbar \Phi (x_{\max})\le 0$ by the maximum principle. Since $\Theta^{i\bar j}$ is positive definite, we calculate at $x_{\max}$,

\bea
0 
& \ge& \nonumber \Theta^{i\bar j} \Phi_{\bar j i} \\ 
& = & \nonumber \frac{n\varepsilon}{n+1} (-\psi_{s,k})^{-\frac{1}{n+1}} \Theta^{i\bar j} (\psi_{s,k})_{i\bar j} + \frac{\varepsilon n}{(n+1)^2} (-\psi_{s,k})^{-\frac{n+2}{n+1}} \Theta^{i\bar j} (\psi_{s,k})_i (\psi_{s,k})_{\bar j} \\
 && \nonumber \quad  - \beta + \tr_{G} \omega_h + \tr_{G}\chi -  \frac{\epsilon'}{n-1} ( \tr_{G} \omega \cdot \tr_{\omega} \omega_{{\mathbb C}^n} - \tr_{G} \omega_{{\mathbb C}^n} ) \\
 & \ge & \nonumber \frac{n\varepsilon}{n+1} (-\psi_{s,k})^{-\frac{1}{n+1}} \Theta^{i\bar j} (\psi_{s,k})_{i\bar j}  - \beta + \tr_{G} [ \omega_h - \frac{2\epsilon'}{n-1} (\tr_\omega \omega_h) \omega ] \\
 & \ge & \nonumber \frac{n^2\varepsilon}{n+1} (-\psi_{s,k})^{-\frac{1}{n+1}} (\det \Theta^{i\bar j})^{1/n} [\det ( \psi_{s,k})_{i\bar j}]^{1/n} - \beta\\
 &\ge& \frac{n^2\varepsilon}{n+1} 
 (-\psi_{s,k})^{-\frac{1}{n+1}}\frac{\gamma^{1/n}}{f}  \frac{(-u_s)^{1/n}e^F}{A_{s,k}^{1/n}} - \beta \nonumber\\
 & \ge & \nonumber \frac{n^2\varepsilon}{n+1} \gamma^{1/n}
 (-\psi_{s,k})^{-\frac{1}{n+1}}\frac{(-u_s)^{1/n}}{A_{s,k}^{1/n}} - \beta.
 \eea

 Here the first equality follows from the definition of $\Theta^{i\bar j}$ and part (a) of Lemma \ref{lm:theta}
 $$\Theta^{i\bar j}\varphi_{i\bar j} \leq \beta - G^{i\bar j}(g_h)_{\bar j i} -  G^{i\bar j}\chi_{\bar j i} .$$ 
 The third line follows from the choice of $\epsilon'$ in \neweqref{eqn:epsilon small}.
In the fourth line, we applied the standard arithmetic-geometric inequality. The fifth line follows from part (b) of Lemma \ref{lm:theta} and the definition of $\psi_{s,k}$. The last line follows from the equation \neweqref{eqn:main}. By the choice of $\varepsilon$, this implies that $\Psi(x_0)\le 0$. Hence $\sup_\Omega \Psi\le 0$. Q.E.D.

\medskip

\bigskip

Along the same spirit in \cite{GP}, as long as we establish the comparison between $u_s$ and $\psi_{s,k}$ as in Lemma \ref{lm:comparison}, we could derive the $L^\infty$ estimate, without refering to the differential equations satisfied by $u_s$ and $\psi_{s,k}$. We cite the lemma here without repeating the proof from \cite[Lemma 2]{GP}.

\begin{lemma}
\label{general}
Assume the functions $u_s= \vp - \vp(x_0) + q(z) - s$ satisfies $u_s>0$ on $\p\Omega$ and $u_s(x_0) <0$. Let $\Omega_s$, $A_s$, and $A_{s,k}$ be the corresponding notions defined above. 

Assume that
\bea \label{eqn:lmgeneral}
-u_s\leq C(n,\gamma,\beta)\,A_{s,k}^{1\over n+1}(-\psi_{s,k})^{n\over n+1}\ {\rm on}\ \bar\Omega
\eea
for some constant $C(n,\gamma,\beta)$, where $\psi_{s,k}$ are plurisubharmonic functions on $\Omega$ such that $\int_\Omega(i\p\bar\p\psi_{s,k})^n=1$ and $\psi_{s,k}=0$ on $\p\Omega$. Then
for any $p>n$, we have
\bea
-\varphi(x_0)\leq C(n,\omega, \gamma,\beta, p, \|\varphi\|_{L^1(\Omega,\o^n)}).
\eea
\end{lemma}

\bigskip

Now, Theorem \ref{local} follows from Lemma \ref{lm:comparison} and Lemma \ref{general}, except that the a priori $L^\infty$ bound of $\varphi$ may rely on the $L^1$ norm of $\varphi$. To remove this dependence, we cite the following lemma from \cite[Lemma 8]{GP}.

\begin{lemma}\label{lemma 4}
Let $\varphi\in C^2(X)$ so that $\lambda[\o^{-1}\cdot (\o + i\partial \bar \partial \varphi)]\in \Gamma_1$ with $\sup_X \varphi =0$, then
$$\int_X (-\varphi) \omega^n \le C(n,\omega).$$
\end{lemma}

\noindent \textbf{Remark.} We stress again that our definition of relative endomorphism $h_\varphi := \o ^{-1} \tilde \o$ is different from those in \cite{GPT,GPTa,GPTW,GPTW1,GPS,GP} due to a different form of unknown metric \neweqref{eqn:unknownmetric}. So the condition $\lambda \in \Gamma$ means $\lambda[\o^{-1}\cdot \tilde \o] \in \Gamma $, instead of $\lambda[\o^{-1}\cdot (\o + i\partial \bar \partial \varphi)]\in \Gamma$. 

Since our assumption for $\varphi$ is that $\lambda[\o^{-1}\cdot \tilde \o]$ lies in $\Gamma \subset \Gamma_1$, i.e. $\tr_\o \tilde \o  \ge 0 $. This implies that $\Delta_{\omega} \varphi \ge - \tr_{\omega} \omega_h \ge -C'$ by \neweqref{eqn:unknownmetric}, for some positive constant $C'$. Therefore $\tr_\o(\o+ i\partial \bar \partial (\frac{n\varphi}{C'})) \ge 0$. In other words, $\lambda[\o^{-1}\cdot (\o + i\partial \bar \partial (\frac{n\varphi}{C'}))]$ is in $\Gamma_1$. Lemma \ref{lemma 4} then yields the desired $L^1$ bound of $\varphi$ given our normalization $\sup_X \varphi = 0$. The proof of Theorem \ref{thm:local} is complete.

\bigskip

\noindent Department of Mathematics, Columbia University, New York, NY 10027 USA

\noindent nklemyatin@math.columbia.edu

\medskip

\noindent Department of Mathematics, Columbia University, New York, NY 10027 USA

\noindent sliang@math.columbia.edu

\medskip

\noindent Department of Mathematics, Columbia University, New York, NY 10027 USA

\noindent wangchuwen@math.columbia.edu


\medskip

\begin{thebibliography}{99}

{\footnotesize












\bibitem{CKNS} L. Caffarelli, J. J. Kohn, L. Nirenberg, and J. Spruck, ``The Dirichlet problem for nonlinear second-order elliptic equations. II. Complex Monge-Amp\`ere, and uniformly elliptic, equations''. Comm. Pure Appl. Math. 38 (1985), no. 2, 209 - 252.

\bibitem{CNS} L. Caffarelli, L. Nirenberg, and J. Spruck, ``The Dirichlet problem for nonlinear second order elliptic equations,
III: functions of the eigenvalues of the Hessian", Acta Math. 155 no. 3-4 (1985) 261-301.


\bibitem{CC2} X.X. Chen and J.R. Cheng, ``The $L^\infty$ estimates for parabolic complex Monge-Amp\`ere and Hessian equations'', arXiv:2201.13339.

\bibitem{CX} J.R. Cheng and Y.L. Xu, ``Regularization Of $ m $-subharmonic Functions And H\" Older Continuity", arXiv:2208.14539.










\bibitem{DK2} S. Dinew and S. Kolodziej, ``A priori estimates for complex Hessian equations",
Anal. PDE 7 no 1 (2013) 227-244.








\bibitem{GPT} B. Guo, D.H. Phong, and F. Tong, ``On $L^\infty$ estimates for complex Monge-Amp\`ere equations'',   arXiv:2106.02224


\bibitem{GPTa} B. Guo, D.H. Phong, and F. Tong, ``Stability estimates for the complex Monge-Amp\`re and Hessian equations'',  arXiv:2106.03913

\bibitem{GPTW} B. Guo, D.H. Phong, F. Tong, and C. Wang, ``On $L^\infty$ estimates for Monge-Amp\`ere and Hessian equations on nef classes'',  arXiv:2111.14186

\bibitem{GPTW1}B. Guo, D.H. Phong, F. Tong, and C. Wang,  ``On the modulus of continuity of solutions to complex Monge-Amp\`ere equations'',  arXiv:2112.02354

\bibitem{GPS} B. Guo, D.H. Phong, and J. Sturm, ``Green's functions and complex Monge-Amp\`ere equations", arXiv:2202.04715.

\bibitem{GP} B. Guo and D.H. Phong, ``On $L^\infty$ estimates for fully nonlinear partial differential equations'', preprint arXiv:2204.12549

\bibitem{GPa} B. Guo and D.H. Phong, ``Uniform entropy and energy bounds for fully non-linear equations'', preprint arXiv:2207.08983

\bibitem{GPb} B. Guo and D.H. Phong, ``Auxiliary Monge-Ampere equations in geometric analysis'', preprint arXiv:2210.13308

\bibitem{GPSS} B. Guo, D.H. Phong, J. Song, and J. Sturm, ``Diameter estimates in K\"ahler geometry'', preprint arXiv:2209.09428

\bibitem{GS} B. Guo and J. Song,  ``Local noncollapsing for complex Monge-Amp\`ere equations'', to appear in {J. Reine Angew. Math.} arXiv:2201.02930

\bibitem{HL} F. R. Harvey and H. B. Lawson,  ``Dirichlet duality and the nonlinear Dirichlet problem''. Comm. Pure Appl. Math. 62 (2009), no. 3, 396 - 443. 

\bibitem{HLnew}  F. R. Harvey and H. B. Lawson,  ``Dirichlet duality and the nonlinear Dirichlet problem on Riemannian manifolds'', J. Differential Geom. 88 (2011), no. 3, 395 - 482.

\bibitem{HLnew2} Harvey, F. R. and Lawson, H. B.,  {\em Geometric plurisubharmonicity and convexity: an introduction}, Adv. Math. {\bf 230} (2012), no. 4-6, 2428--2456.

\bibitem{HL3} F.R. Harvey and H.B. Lawson, ``Determinant majorization and the work of Guo-Phong-Tong and Abja-Olive", arXiv: 2207.01729.

\bibitem{K} S. Kolodziej, ``The complex Monge-Amp\`ere equation", Acta Math. 180 (1998) 69-117.






\bibitem{S}M. Sroka, ``Sharp uniform bound for the quaternionic Monge-Amp\'ere equation on hyperhermitian manifolds", arXiv preprint arXiv:2211.00959.

\bibitem{STW} G. Szekelyhidi, V. Tosatti, and B. Weinkove,
``Gauduchon metrics with prescribed volume form",
Acta Math. 219 (2017) no. 1, 181-211.



\bibitem{TW1} V. Tosatti, and B. Weinkove, ``The Monge-Amp\`ere equation for $(n-1)$-plurisubharmonic functions on a compact K\"ahler manifold'', J. Amer. Math. Soc. 30 (2017), no.2, 311-346.


\bibitem{TW} N.S. Trudinger and X.J. Wang, ``Hessian equations II", Ann. of Math. 150 (1999) 579-604. 




\bibitem{Y} S.T. Yau, ``On the Ricci curvature of a compact K\"ahler manifold and the complex Monge-Amp\`ere equation. I'', Comm. Pure Appl. Math. 31 (1978) 339-411.



}

\end{thebibliography}
\end{document}